\documentclass[a4paper,12pt]{amsart}

\usepackage{amsmath,amsfonts,amssymb,amsthm}

\usepackage{subfig}
\usepackage{graphicx}
\usepackage{color}
\usepackage[all]{xy}
\usepackage{url}
\usepackage{setspace}
\usepackage{booktabs}
\usepackage{longtable}
\usepackage{latexsym}
\usepackage{comment}
\usepackage[usenames,dvipsnames]{xcolor}

\usepackage{tikz}
\usetikzlibrary{cd}
\usetikzlibrary{patterns}
\usetikzlibrary{shapes}
\usetikzlibrary{decorations.fractals}

\usepackage{hyperref}
\hypersetup{
    colorlinks = true,
    linkcolor  = {red!50!black},
    citecolor  = {blue!50!black},
    urlcolor   = {blue!80!black}
}

\newtheorem*{introthm}{Theorem}
\newtheorem{theorem}{Theorem}

\theoremstyle{definition}
\newtheorem{definition}[theorem]{Definition}

\definecolor{cKlaus}{rgb}{0.1,0.55,0.03}
\definecolor{cKlausOK}{rgb}{0.6,0.10,0.33}
\definecolor{intOrange}{rgb}{1.0,.310,.0}

\newcommand{\ktrash}[1]{}

\newcommand{\dtrash}[1]{}
\usepackage[textsize=tiny]{todonotes}

\DeclareMathOperator{\tail}{tail}

\DeclareMathOperator{\orb}{orb}

\DeclareMathOperator{\Spec}{Spec}

\DeclareMathOperator{\conv}{conv}

\newcommand{\kprojlimm}[1]{\underset{#1}{\varprojlim}}

\newcommand{\coloneqq}{\mathrel{\mathop:}=}

\newcommand{\PP}{\mathbb P}
\newcommand{\A}{\mathbb A}

\newcommand{\F}{\mathbb F}

\newcommand{\N}{\mathbb N}
\newcommand{\R}{\mathbb R}
\newcommand{\T}{\mathbb T}
\newcommand{\V}{\mathbb V}
\newcommand{\Z}{\mathbb Z}

\newcommand{\CN}{{\mathcal N}}
\newcommand{\CL}{{\mathcal L}}
\newcommand{\CO}{{\mathcal O}}

\newcommand{\CS}{{\mathcal S}}

\newcommand{\CU}{{\mathcal U}}

\newcommand{\wt}{\widetilde}
\newcommand{\pairing}[2]{\langle#1,#2\rangle}
\newcommand{\til}[1]{\widetilde{#1}}

\newcommand{\ifff}{\,\Longleftrightarrow\,}

\newcommand{\gExt}{\operatorname{Ext}}

\newcommand{\gHom}{\operatorname{Hom}}

\newcommand{\gH}{\operatorname{H}}
\newcommand{\tH}{\wt{\gH}}

\newcommand{\kst}{\,|\;}

\newcommand{\kss}{\scriptscriptstyle}
\newcommand{\kbb}{{\kss \bullet}}
\newcommand{\ko}{\overline}

\newcommand{\Pol}{\operatorname{Pol}}
\newcommand{\dual}{^{\scriptscriptstyle\vee}}

\newcommand{\kk}{k}

\newcommand{\normal}{\CN}
\newcommand{\toric}{\T\V}  
\newcommand{\Pic}{\operatorname{Pic}}
\newcommand{\Cl}{\operatorname{Cl}}
\newcommand{\Db}{\operatorname{D}^b} 
\newcommand{\Coh}{\operatorname{Coh}}

\definecolor{skin}{HTML}{FFECC9}
\definecolor{pumpkin}{HTML}{FEDFA9}
\definecolor{piggy}{HTML}{FFB99D}
\definecolor{fiolet}{HTML}{CD8F9C}
\definecolor{granat}{HTML}{677081}
\definecolor{ciemnyblekit}{HTML}{91A1B8}

\definecolor{oliwkowy}{HTML}{627037}
\definecolor{ciemnazielen}{HTML}{394D2E}
\definecolor{ciemnyfiolet}{HTML}{424444}
\definecolor{mocnyfiolet}{HTML}{717299}
\definecolor{jasnyfiolet}{HTML}{B0ABCC}
\definecolor{bladyfiolet}{HTML}{C9C7DB}

\definecolor{lightblue}{RGB}{135,206,250}
\definecolor{darkblue}{RGB}{0,0,160}
\definecolor{darkgreen}{RGB}{0,100,0}

\newcommand{\bild}[2]{\begin{tikzpicture}[scale=#1]#2\end{tikzpicture}}

\newcommand{\PolPlusTail}{\Pol^+_{\,\delta}}
\newcommand{\PolTail}{\Pol_{\,\delta}}
\newcommand{\PolPlusSigma}{\PolPlusdeltaSigma}
\newcommand{\PolPlusdeltaSigma}{\Pol^+_\delta(\Sigma)}
\newcommand{\PolPlusnullSigma}{\Pol^+_0(\Sigma)}

\newcommand{\PoldeltaSigma}{\Pol_\delta(\Sigma)}
\newcommand{\CLL}{\CO}

\newcommand{\vp}{v^+_\sigma}
\newcommand{\vm}{v^-_\sigma}
\newcommand{\vv}{v_\sigma}
\newcommand{\vvPrime}{v_{\sigma'}}
\newcommand{\vf}{v^\kbb_\sigma}
\newcommand{\vpm}{v^\Delta_\sigma}

\begin{document}

\title[Displaying the cohomology of toric line bundles]
{Displaying the cohomology of toric line bundles}

\author[K.~Altmann]{Klaus Altmann%
}
\address{Institut f\"ur Mathematik,
FU Berlin,
Arnimallee 3,
14195 Berlin,
Germany}
\email{altmann@math.fu-berlin.de}
\author[D.~Ploog]{David Ploog}
\address{Fachbereich Mathematik,
Welfengarten 1, 30167 Hannover, Germany}
\email{ploog@math.uni-hannover.de}

\begin{abstract}
There is a standard method to calculate the cohomology of torus-invariant
sheaves $\CL$ on a toric variety via the simplicial cohomology of  
associated subsets $V(\CL)$ of the space $N_\R$ of 1-parameter subgroups 
of the torus. 
For a line bundle $\CL$
represented by a formal difference $\Delta^+-\Delta^-$
of polyhedra in the character space $M_\R$,
\cite{immaculate} contains a simpler formula for the cohomology of $\CL$, 
replacing $V(\CL)$ by the set-theoretic difference $\Delta^- \setminus \Delta^+$.
Here, we provide a short and direct proof of this formula.
\end{abstract}

\maketitle

\section{Introduction}
\label{intro}

\subsection{Toric varieties}\label{tVar}
Let $\kk$ be an algebraically closed field. We consider the algebraic torus
$T=\Spec\kk[M]$ where $M$ is a free abelian group of finite rank
and $\kk[M]$ denotes the semigroup ring of $M$.
Then $M$ can be recovered from $T$ as its character group. 
The dual $N = \gHom_{\Z}(M,\Z)$
is the group of 1-parameter subgroups in $T$.
By definition, we have a natural perfect pairing
$$
\langle-,-\rangle \colon M\times N \to \Z.
$$
The theory of toric varieties deals with partial compactifications
$j \colon T\hookrightarrow X$ such that the group law of $T$ extends to an algebraic action of
$T$ on the $\kk$-variety $X$. By a \emph{cone} in 
$N_\R \coloneqq N\otimes_Z\R$, we mean a 
subset $\sigma\subseteq N_\R$ which is finitely generated and convex
($\sigma = \sum_{i=1}^r\R_{\geq 0} \cdot v_i$) and
rational ($v_i\in N$).
Each affine toric variety is obtained from a 
{\em pointed} cone $\sigma\subseteq N_\R$, i.e.\ $\sigma\cap(-\sigma)=\{0\}$, via
$$
\toric(\sigma,N) = \toric(\sigma) \coloneqq \Spec\kk[\sigma\dual\cap M]
\hspace{0.8em}\mbox{with}\hspace{0.8em}
\sigma\dual \coloneqq \{m\in M_\R
\kst \langle m,\sigma\rangle\geq 0\}.
$$
Each face $\tau\leq\sigma$ provides a natural open
embedding $\toric(\tau)\hookrightarrow\toric(\sigma)$. 
In particular, $\{0\}\leq\sigma$ corresponds
to the open subset $T = \toric(0) \hookrightarrow \toric(\sigma)$.
Gluing affine toric varieties along such open subsets leads to general (non-affine) toric varieties
determined by a polyhedral fan $\Sigma$ in $N_\R$:
$$
X = \toric(\Sigma) \coloneqq \lim_{\substack{\longrightarrow \\ \sigma\in\Sigma}} \toric(\sigma) .
$$
A {\em fan} $\Sigma$ is a finite set of pointed cones in $N_\R$ 
closed under taking faces and such that 
$\sigma\cap\sigma'\leq\sigma, \sigma'$ for any two cones
$\sigma,\sigma'\in\Sigma$.
Then the affine varieties $\toric(\sigma)$ turn into affine, open charts
$U_\sigma\subseteq X=\toric(\Sigma)$.
\\[1ex]
The main advantage of toric varieties is that the $T$-action provides
a fine $M$-grading on all algebraic structures functorially associated to
$\toric(\Sigma)$, 
allowing a combinatorial description. 
Examples are the cohomology groups of $T$-invariant sheaves or the modules 
describing the infinitesimal deformations and obstructions of toric varieties.
The latest general reference for the theory of toric varieties is
the book \cite{CoxBook}.

\subsection{Cohomology of $T$-invariant Weil divisors}\label{cohomTinvWeil}
Let us assume that the {\em support} $|\Sigma| \coloneqq \bigcup_{\sigma\in\Sigma}\sigma$ 
of a given fan $\Sigma$
is a full-dimensional, not necessarily pointed cone, 
e.g.\ $|\Sigma| = N_\R$.
We will identify the rays, i.e.\ the one-dimensional cones $\rho\in\Sigma(1)$,
with their primitive generators $\rho\in N$. 
They give rise to the one-codimensional $T$-orbits 
$\orb(\rho)  \coloneqq \Spec\kk[\rho^\bot\cap M] \subseteq
 \toric(\rho)\subseteq X$,
and their closures
$$
D_\rho \coloneqq \ko{\orb(\rho)}
$$
are precisely the $T$-invariant (or ``toric'') prime divisors. 
Hence, an arbitrary toric Weil divisor has the form
$$
D=\sum_{\rho\in\Sigma(1)}\lambda_\rho \,D_\rho
\hspace{1em}(\lambda_\rho\in\Z).
$$
One of the salient features of toric varieties is that these special divisors $D_\rho$
generate the full class group. This is reflected by the famous exact sequence
$$
0 \to M \to \Z^{\Sigma(1)} \to \Cl(X)\to 0
$$
whose first map 
is the dual of the natural map $\Z^{\Sigma(1)}\to N$, $e_\rho\mapsto\rho\in N$,
sending standard basis vectors to primitive ray generators.
During the last five decades, it has been one of the 
basic results in toric geometry to express the 
$M$-eigenspaces of the sheaf
cohomology of $\CO_X(D)$ on $X=\toric(\Sigma)$ as the singular cohomology of certain
subsets $V_{D,m}$ of the vector space $N_\R$, where $m\in M$. In 
\cite[Prop.~2.6]{demazure},
\cite[I.3]{KKMS},
\cite[\S7]{danilov},
\cite[(2.2)]{OdaBook},
\cite[(3.5)]{fultonToric}, or in
\cite[(9.1)]{CoxBook}
one defines, for all characters $m\in M$ these sets as
$$
V_{D,m} \coloneqq \bigcup_{\sigma\in\Sigma}
\conv\{\rho\in\sigma(1)\kst \langle m,\rho\rangle < -\lambda_\rho\}
\subseteq N_\R\setminus\{0\}.
$$
The basic formula allured to above is, with coefficients $\kk$ on
the right-hand side
$$
\gH^i\!\big(X,\CO_X(D)\big)(m)=
\tH^{i-1}(V_{D,m}).
$$
Recall that the $(-1)$-st reduced cohomology is 
$\tH^{-1}(V)=\left\{\begin{array}{ll}
0 & \mbox{if } V\neq\emptyset\\
\kk & \mbox{if } V=\emptyset
\end{array}\right.$.

\subsection{Cohomology of $T$-invariant Cartier divisors}\label{cohomTinvCart}
While the previous formula is very useful, its main ingredient
is the rather technically defined subset $V_{D,m} \subseteq N_\R$. 
In \cite[III.3]{immaculate},
a Fourier transformation argument was used to replace these subsets
by much more natural subsets of the dual vector space
$M_\R$, at least for Cartier divisors in quasi-projective toric varieties.
This makes it possible to literally display the cohomology of the line bundles 
representing the divisors:

In this setting, a nef Cartier divisor $D$ is represented by a convex 
lattice polyhedron $\Delta$ in $M_\R$,
see Subsection~(\ref{polToCartier}).
The lattice points of the polyhedron
correspond to global sections of $\CO_X(D)$,
yielding $k^{\Delta\cap M} = \gH^0(X,\CO_X(D))$.

An arbitrary Cartier divisor is a difference $D = D^+-D^-$ of two nef
Cartier divisors, and then these two parts 
can be represented by convex lattice polyhedra
$
\Delta^+, \Delta^-\subseteq M_\R,
$
respectively.
The Minkowski sums among those polytopes correspond to the sums
of the associated divisors, i.e.\ to the tensor products of their
sheaves of sections.
The divisor $D$ itself corresponds to the formal difference
$\Delta^+-\Delta^-$ in the Grothendieck group of the semigroup
of polyhedra with Minkowski addition, 
see Section \ref{cartPol} for more details. Now, the result of
Theorem III.6 in \cite{immaculate}
is

\begin{introthm}
\label{th-cohomByDiff}
$\;\gH^i\!\big(X,\CO_X(D)\big)(m)=
\tH^{i-1}\!\big(\Delta^-\setminus(\Delta^+-m)\big)$.
\end{introthm}

The goal of the present paper is to provide a direct, straight-forward proof of this result 
(subsequently called Theorem~\ref{th-cohomByDiffB}), avoiding the application
of the traditional formula and of the Fourier transformation argument linking
$V_{D,m}\subset N_\R$ to $M_\R$. This might become especially
useful when trying to generalize Theorem~\ref{th-cohomByDiffB} to the situation
of Okounkov bodies or to $T$-varieties of higher complexity. In both
situations, the ``$N$-side'' of the story does not exist or is at least not
easily accessible.

\subsection{An easy example}\label{easyEx}
The first Hirzebruch surface
$\F_1=\PP\big(\CO_{\PP^1}\oplus\CO_{\PP^1}(1)\big)$ is obtained by blowing
up $\PP^2$ in one (torus-invariant) point. The fan of $\F_1$ looks like
$$
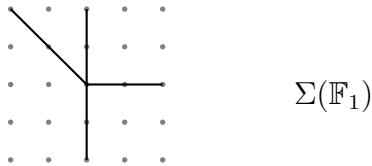
\bild{0.5}{
\foreach \x in {-2,-1,...,2} \foreach \y in {-2,-1,...,2} {
  \fill[thick, color=gray] (\x,\y) circle (2pt); }
\draw[thick,  color=black]
  (0,2) -- (0,-2) (0,0) -- (2,0) (0,0) -- (-2,2);
\draw[thick,  color=black]  (6.5,-0.3) node{$\Sigma(\F_1)$};
}
$$
The nef cone in $\Cl(\F_1)\cong\Z^2$
is freely generated by the line bundles associated to the two
polytopes
$$
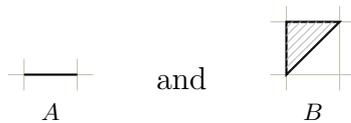
\newcommand{\sizeX}{0.7}
\bild{\sizeX}{
\draw[color=oliwkowy!40] (-0.3,-0.3) grid (1.3,0.3);
\draw[thick,  color=black]
  (0,0) -- (1,0);
\draw[thick,  color=black]  (0.5,-0.7) node{\scriptsize $A$};
}
\hspace{2em}\raisebox{3ex}{\mbox{and}}\hspace{2em}
\bild{\sizeX}{
\draw[color=oliwkowy!40] (-0.3,-0.3) grid (1.3,1.3);
\draw[thick,  color=black]
  (0,0) -- (0,1) -- (1,1) -- cycle;
\draw[thick,  color=black]  (0.5,-0.7) node{\scriptsize $B$};
\fill[pattern color=gray!50, pattern=north east lines]
  (0,0) -- (0,1) -- (1,1) -- cycle;
}
$$
See Subsection~(\ref{polToCartier}) for details
on this correspondence, and Section~\ref{sub:fes} for the actual
line bundles.
Ample line bundles arise from tensor products of both generators;
the two simplest examples are represented by the Minkowski sums
$$
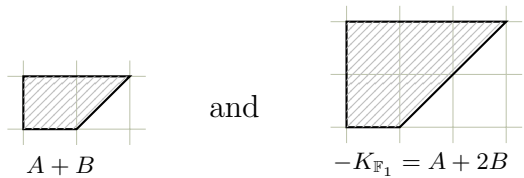
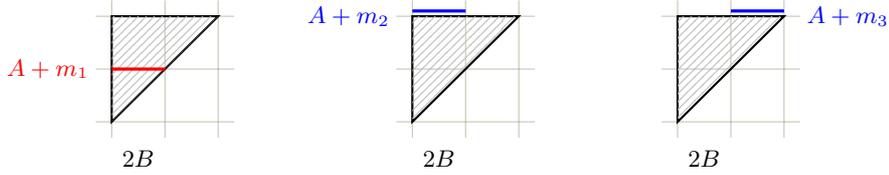
\newcommand{\sizeX}{0.7}
\bild{\sizeX}{
\draw[color=oliwkowy!40] (-0.3,-0.3) grid (2.3,1.3);
\draw[thick,  color=black]
  (0,0) -- (1,0) -- (2,1) -- (0,1) -- cycle;
\draw[thick,  color=black]  (0.7,-0.7) node{\scriptsize $A+B$};
\fill[pattern color=gray!50, pattern=north east lines]
  (0,0) -- (1,0) -- (2,1) -- (0,1) -- cycle;
}
\hspace{2em}\raisebox{5ex}{\mbox{and}}\hspace{2em}
\bild{\sizeX}{
\draw[color=oliwkowy!40] (-0.3,-0.3) grid (3.3,2.3);
\draw[thick,  color=black]
  (0,0) -- (1,0) -- (3,2) -- (0,2) -- cycle;
\draw[thick,  color=black]  (1.4,-0.7) node{\scriptsize $-K_{\F_1}=A+2B$};
\fill[pattern color=gray!50, pattern=north east lines]
  (0,0) -- (1,0) -- (3,2) -- (0,2) -- cycle;
}
$$
Now, the pictures
$$
\newcommand{\sizeX}{0.7}
\bild{\sizeX}{
\draw[color=oliwkowy!40] (-0.3,-0.3) grid (2.3,2.3);
\draw[thick,  color=black]
  (0,0) -- (0,2) -- (2,2) -- cycle;
\draw[thick,  color=black]  (0.5,-0.7) node{\scriptsize $2B$};
\fill[pattern color=gray!50, pattern=north east lines]
  (0,0) -- (0,2) -- (2,2) -- cycle;
\draw[very thick,  color=red] (0,1) -- (1,1);
\draw[thick,  color=red]  (-1.2,1.0) node{\scriptsize $A+m_1$};
}
\hspace{2em}
\bild{\sizeX}{
\draw[color=oliwkowy!40] (-0.3,-0.3) grid (2.3,2.3);
\draw[thick,  color=black]
  (0,0) -- (0,2) -- (2,2) -- cycle;
\draw[thick,  color=black]  (0.5,-0.7) node{\scriptsize $2B$};
\fill[pattern color=gray!50, pattern=north east lines]
  (0,0) -- (0,2) -- (2,2) -- cycle;
\draw[very thick,  color=blue] (0,2.1) -- (1,2.1);
\draw[thick,  color=blue]  (-1.2,2.0) node{\scriptsize $A+m_2$};
}
\hspace{4em}
\bild{\sizeX}{
\draw[color=oliwkowy!40] (-0.3,-0.3) grid (2.3,2.3);
\draw[thick,  color=black]
  (0,0) -- (0,2) -- (2,2) -- cycle;
\draw[thick,  color=black]  (0.5,-0.7) node{\scriptsize $2B$};
\fill[pattern color=gray!50, pattern=north east lines]
  (0,0) -- (0,2) -- (2,2) -- cycle;
\draw[very thick,  color=blue] (1,2.1) -- (2,2.1);
\draw[thick,  color=blue]  (3.2,2.0) node{\scriptsize $A+m_3$};
}
$$
show that there are exactly three degrees $m_i\in M=\Z^2$ 
such that the translates $A+m_i$ are contained in $2B$,
i.e.\ that $(A+m_i)\setminus (2B)=\emptyset$ -- but only the shift by
$m_1$ divides $(2B)\setminus(A+m_1)$ into two connected components.
Hence
$\gH^0\!\big(\F_1,\,\CO_{\F_1}(2B-A)\big)$ is three-dimensional,
supported in the degrees $m_1,m_2,m_3\in M$, and
$$
\gH^1\!\big(\F_1,\,\CO_{\F_1}(A-2B)\big) =
\tH^0\!\big(2B\setminus (A+m_1)\big)
$$
is one-dimensional,
supported in the single degree $-m_1\in M$.
Similarily, one can spot the unique
integral shift of $A$ 
into the interior of $4B$ leading to a one-dimensional
$\gH^2\big(\F_1,\,\CO_{\F_1}(A-4B)\big)$ supported in $-m_4\in M$:
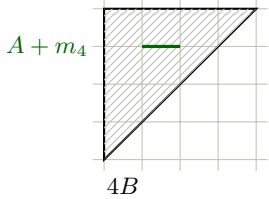
$$
\newcommand{\sizeX}{0.5}
\bild{\sizeX}{
\draw[color=oliwkowy!40] (-0.3,-0.3) grid (4.3,4.3);
\draw[thick,  color=black]
  (0,0) -- (0,4) -- (4,4) -- cycle;
\draw[thick,  color=black]  (0.5,-0.7) node{\scriptsize $4B$};
\fill[pattern color=gray!50, pattern=north east lines]
  (0,0) -- (0,4) -- (4,4) -- cycle;
\draw[very thick,  color=darkgreen] (1,3) -- (2,3);
\draw[thick,  color=darkgreen]  (-1.5,3.0) node{\scriptsize $A+m_4$};
}
$$
Since $-K_{\F_1}=A+2B$, one can easily visualize Serre duality.
For instance, the one-dimensional
$\gH^1\!\big(\F_1,\,\CO_{\F_1}(A-2B)\big)$ leads to the
one-dimensional
$\gH^1\!\big(\F_1,\,\CO_{\F_1}(-2A)\big)^{\!\vee}$, and the latter can be
seen from Theorem~\ref{th-cohomByDiffB} since the polytope $2A$ is an
interval of length 2 and thus has a unique interior lattice point.

\subsection{An even easier example with non-trivial tail cone}\label{easyerEx}
The blowing up $\til{\A}^2$
of the origin in $\A^2$ is represented by the natural map of fans
$$
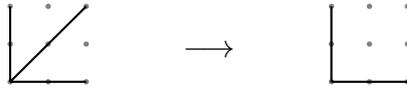
\newcommand{\sizeX}{0.5}
\bild{\sizeX}{
\foreach \x in {0,...,2} \foreach \y in {0,...,2} {
  \fill[thick, color=gray] (\x,\y) circle (2pt); }
\draw[thick,  color=black]
  (0,2) -- (0,0) (0,0) -- (2,0) (0,0) -- (2,2);
}
\hspace{3em}\raisebox{2ex}{$\longrightarrow$}\hspace{3em}
\bild{\sizeX}{
\foreach \x in {0,...,2} \foreach \y in {0,...,2} {
  \fill[thick, color=gray] (\x,\y) circle (2pt); }
\draw[thick,  color=black]
  (0,0) -- (2,0) (0,0) -- (0,2);
}
$$
The negative of the exceptional divisor $E\subseteq\til{\A}^2$ 
is ample, and the associated sheaf $\CO(-E)$ is represented by the polyhedron
(see Subsection~(\ref{polNormFan}) for tail cones)
$$
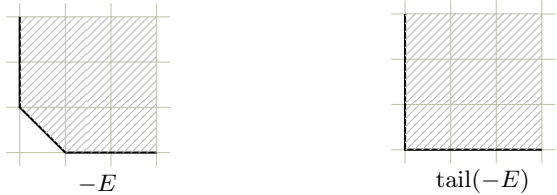
\newcommand{\sizeX}{0.6}
\bild{\sizeX}{
\draw[color=oliwkowy!40] (-0.3,-0.3) grid (3.3,3.3);
\draw[thick,  color=black]
  (0,3) -- (0,1) -- (1,0) -- (3,0);
\draw[thick,  color=black]  (1.7,-0.7) node{\scriptsize $-E$};
\fill[pattern color=gray!50, pattern=north east lines]
  (0,3) -- (0,1) -- (1,0) -- (3,0) -- (3,3) -- cycle;
}
\hspace{7em}
\bild{\sizeX}{
\draw[color=oliwkowy!40] (-0.3,-0.3) grid (3.3,3.3);
\draw[thick,  color=black]
  (0,3) -- (0,0) -- (3,0);
\draw[thick,  color=black]  (1.7,-0.7) node{\scriptsize $\tail(-E)$};
\fill[pattern color=gray!50, pattern=north east lines]
  (0,3) -- (0,0) -- (3,0) -- (3,3) -- cycle;
}
$$
Now, Theorem~\ref{th-cohomByDiffB} and this picture
$$
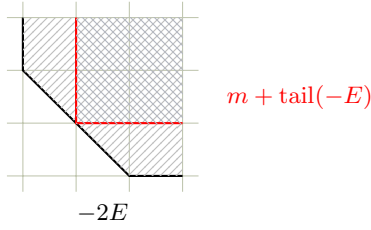
\newcommand{\sizeX}{0.7}
\bild{\sizeX}{
\draw[color=oliwkowy!40] (-0.3,-0.3) grid (3.3,3.3);
\draw[thick,  color=black]
  (0,3) -- (0,2) -- (2,0) -- (3,0);
\draw[thick,  color=black]  (1.5,-0.7) node{\scriptsize $-2E$};
\fill[pattern color=gray!50, pattern=north east lines]
  (0,3) -- (0,2) -- (2,0) -- (3,0) -- (3,3) -- cycle;
\draw[thick,  color=red]
  (1,3) -- (1,1) -- (3,1);
\fill[pattern color=granat!50, pattern=north west lines]
  (1,3) -- (1,1) -- (3,1) -- (3,3) -- cycle;
\draw[thick,  color=red]  (5.2,1.5) node{\scriptsize $m+\tail(-E)$};
}
$$
show that $\gH^1\!\big(\til{\A}^2,\, \CO_{\til{\A}^2}(0E-(-2E))\big)
= \gH^1\!\big(\til{\A}^2,\, \CO_{\til{\A}^2}(2E)\big)$ is one-dimensional.

\section{Cartier divisors and polyhedra}
\label{cartPol}

\subsection{Lattice polytopes with prescribed normal fan}\label{polNormFan}
If $\Delta\subseteq M_\R$ is a lattice polyhedron, then we define
its \emph{tail cone} (also called its recession cone) as
$$
\tail(\Delta) \coloneqq \{m\in M_\R\kst m+\Delta\subseteq\Delta\}.
$$
Now, we turn the tables. Once and for all, we fix a pointed 
cone $\delta\subseteq M_\R$ and take it as the 
prescribed tail cone of all our polyhedra, i.e.\ we consider the set
$$
\PolPlusTail \coloneqq \{\mbox{lattice polyhedra }\Delta\subseteq M_\R\kst
\tail(\Delta)=\delta\}.
$$
The most important case is $\delta=0$, as $\Pol^+_0$ 
is the set of compact lattice polyhedra. 
In any case, $\PolPlusTail$ is a semigroup under Minkowski addition
with neutral element $\delta$.

\begin{definition}
\label{def-normalFan}
Let $\Delta\subseteq M_\R$ be a lattice polyhedron with
$\tail(\Delta)=\delta$. A fan $\Sigma$ in $N_\R$ is called 
\emph{compatible with $\Delta$} if it is a subdivision of the normal
fan $\normal(\Delta)$, i.e.\ if $|\Sigma|=\delta\dual\subseteq N_\R$,
and if the concave, piecewise linear function
$$
\delta\dual
\to\R,\hspace{1em}
n\mapsto \min \langle\Delta,\,n\rangle
$$
is linear on the cones of $\Sigma$.
\end{definition}

We turn the tables once more:
Fixing a fan $\Sigma$ with support $|\Sigma|=\delta\dual$, the subset
$$
\PolPlusdeltaSigma \coloneqq \{\mbox{lattice polyhedra }\Delta\subseteq M_\R\kst
\Delta \mbox{ is compatible with }\Sigma\}
\subseteq\PolPlusTail
$$
is a finitely generated subsemigroup.
Let us now fix a polyhedron $\Delta\in\PolPlusSigma$,
i.e.\ a compatible pair $(\Delta,\Sigma)$.
{By definition, the function
$\min \langle\Delta,\,\kbb\rangle$
is linear on every cone $\sigma\in\Sigma$. If $\sigma$ is full-dimensional,} this function is realized by a unique element
$$
\vv = \vv^\Delta \in M,
\hspace{0.4em}\text{i.e., for all $n\in\sigma$, we have}\hspace{0.4em}
\min\langle\Delta,n\rangle=\big\langle \vv,\,n\big\rangle.
$$
This immediately implies that
$\vv\in \Delta$ is a vertex 
and $\vv + \sigma\dual=\Delta+\sigma\dual$.
Actually,
$$
\Delta = \conv\{\vv\kst \sigma\in\Sigma(\mbox{\rm top})\}
+ \delta.
$$
The functions                        
$
\,\vf \colon \PolPlusSigma\to M$, $\,\Delta \mapsto \vv^\Delta
$
are obviously additive. However,
unless $\normal(\Delta)=\Sigma$, some of them 
will coincide. 
Finally, if $\sigma\in\Sigma$ is not full-dimensional,
then the lattice points $\vv$ still exist, but
they are no longer uniquely determined.
Instead, one obtains a well-defined function
$\vf \colon \PolPlusSigma\to M/(\sigma^\bot\cap M)$.

\subsection{The Grothendieck groups}\label{grothGrp}
Since we have fixed the tail cone, the semigroup $\PolPlusTail$
is cancellative.
That is, $\PolPlusTail$
can be canonically embedded into the group
$$
\PolTail \coloneqq \{\Delta^+-\Delta^-\kst \Delta^{\pm}\in\PolPlusTail\}
$$
of all formal differences.
Analogously, the subsemigroup
$\PolPlusSigma\subseteq\PolPlusTail$ leads to a
finitely generated subgroup 
$$
\PoldeltaSigma=\{\Delta^+-\Delta^-\kst \Delta^{\pm}\in\PolPlusSigma\}
\subseteq\PolTail.
$$
The additive function $\vf$ on $\PolPlusSigma$
extends to $\PoldeltaSigma$, however with  
$$
\vpm \coloneqq \vp - \vm \coloneqq \vv^{\Delta^+} - \vv^{\Delta^-}
$$
now understood as the virtual vertices of the formal difference
$\Delta^+-\Delta^-$.
Recall that in the quasi-projective case of $\Sigma=\normal(\Delta)$,
the set $\PolPlusSigma$ can also be described as
$$
\PolPlusSigma=\{\Delta'\in\PolPlusTail\kst
\exists\Delta''\in\PolPlusTail: 
\Delta'+\Delta''=N\cdot\Delta\mbox{ for some }N\in\N\}.
$$
In particular, in this case,
the elements of $\PoldeltaSigma$ can be written as 
$\Delta^+-N\cdot\Delta$ by just using the 
``ample polyhedron'' $\Delta$ determining $\Sigma$ and some $N\gg 0$ to 
form the negative part.

\subsection{The nef divisor associated to a polyhedron}
\label{polToCartier}
Assume that $\Delta\subseteq M_\R$ is a polyhedron with tail cone $\delta$ and
that $\Sigma$ is a fan in $N_\R$ compatible with $\Delta$.
Consider again the toric variety $X=\toric(\Sigma)$
with its open torus embedding $j\colon T\hookrightarrow X$. 
Denoting by $\chi^m\in\kk[M]$ the monomial associated to $m\in M$,
we define a{n invertible} 
sheaf $\CO_X(\Delta)$ on $X$ via
$$
\CO_X(\Delta)|_{U_\sigma} \;\coloneqq\; \chi^{\vv}\cdot \CO_{U_\sigma}
\;\subseteq\; j_*\CO_T \;=\;\kk[M]
{ \;=\; k[T] \;\subseteq\; k(T) \;=\;} \kk(X)
$$
where the gluing condition for $\CO_X(\Delta)$ on
$U_\sigma\cap U_{\sigma'} = U_{\sigma\cap\sigma'}$ is obtained from
$$
\textstyle
\vv-\vvPrime\in (\sigma\cap\sigma')^\bot,
\;
\mbox{ i.e. }
\chi^{\vv}/\chi^{\vvPrime}\in
\CO_X(U_\sigma\cap U_{\sigma'})^*.
$$
In particular, we obtain
$$
\gH^0\!\big(U_{\sigma},\,\CO_X(\Delta)\big)
= \chi^{\vv}\cdot \kk[\sigma\dual\cap M]
= \kk[{\vv}+(\sigma\dual\cap M)].
$$
The intersection of these groups exhibits $\Delta\cap M$ as a $\kk$-basis
of $\gH^0\!\big(X,\,\CO_X(\Delta)\big)$.
Besides, the assignment $\Delta\mapsto \CO_X(\Delta)$ 
induces the well-known semigroup isomorphism
$$
\PolPlusSigma\stackrel{\sim}{\to}
\{\mbox{\rm globally generated
$T$-equivariant invertible subsheaves of ${k}(X)$}\}.
$$
In fact, this even becomes a functor when defining morphisms among lattice
polyhedra as injections.
An integral shift of $\Delta\in\PolPlusSigma$ by $m\in M$ on the left-hand
side leads to an isomorphic sheaf whose $T$-action is shifted by the character $m$:
$$
\CO_X(\Delta+m)=\chi^m\cdot\CO_X(\Delta) . 
$$
This isomorphism extends to the respective Grothendieck groups.
Moreover, in the quasi-projective case, 
after modding out integral translations on the left-hand side
and linear equivalence, on the right,
this induces the identification
$$
{\PoldeltaSigma}/M \stackrel{\sim}{\to} \Pic\big(\toric(\Sigma)\big)
$$
where the positive cone $\PolPlusSigma$
of globally generated sheaves
maps to the nef cone within $\Pic\big(\toric(\Sigma)\big)$.
See, e.g., \cite[(3.4)]{fultonToric} 
or \cite[Ch.~6]{CoxBook} for more details.

\section{Comparing the cohomology theories}
\label{compCohom}

\subsection{The comparison theorem}\label{compTheorem}
Let $\Delta^+,\Delta^-\subseteq M_\R$ be two polyhedra with tail cone $\delta$
and let $\Sigma$ be a fan in $N_\R$ compatible with both $\Delta^{\pm}$.
This gives rise to the $T$-equivariant, invertible sheaf on $X=\toric(\Sigma)$
$$
\CL \coloneqq
\CO_X(\Delta^+-\Delta^-)=\CO_X(\Delta^+)\otimes_{\CO_X}\CO_X(\Delta^-)^{-1} .
$$

\begin{theorem}
\label{th-cohomByDiffB}
One has
$\;\gH^i\!\big(X,\CL\big)(m)=
\tH{}^{i-1}\big(\Delta^-\setminus(\Delta^+-m)\big)$
for all $m\in M$, $i\in\N$.
\end{theorem}

The rest of this section contains the proof of this statement.
Note that by our remark $\CO_X(\Delta+m)=\chi^m\cdot\CO_X(\Delta)$
at the end of Subsection~(\ref{polToCartier}), we may 
assume $m=0$. Now, everything follows from the comparison
of two \v{C}ech complexes -- and this will be done in the next subsections.

\subsection{The \v{C}ech complex for $\CO_X(\Delta^+-\Delta^-)$}
\label{cechSpectral}
We refer to \cite[III.4]{hartshorne} for the \v{C}ech complex of a sheaf $\CL$ on a topological space.
On $X=\toric(\Sigma)$ we take the open, affine covering
$$
\CU\coloneqq\{U_\sigma\}_{\sigma\in\Sigma}
$$
which is closed under intersections, i.e.\
$ U_{\sigma}\cap U_{\sigma'}=
U_{\sigma\cap\sigma'}$.
We fix an ordering among the participating
polyhedral cones and for each $p\in\N$, we define
$$
\textstyle
C^p(\CU, \CL) \coloneqq \prod\limits_{i_0<\cdots<i_p}\CL(U_{\sigma_{i_0}\cap\ldots\cap\sigma_{i_p}})
$$
together with the usual differentials $d\colon C^p\to C^{p+1}$. The cohomology of
this complex is denoted by $\check{\gH}^p(\CU,\CL)$, and it maps isomorphically 
onto $\gH^p(X,\CL)$ by \cite[Th.~III.4.5]{hartshorne}, since $\CL$ is coherent.
Now, the intersection
$$
\sigma \coloneqq \sigma_{i_0}\cap\ldots\cap\sigma_{i_p} 
$$
is a cone in $\Sigma$, and the $M$-graded $\kk$-vector space of local 
sections on $U_\sigma$ equals
$$
\gH^0\!\big(U_\sigma,\CL\big)=\chi^{\vp}/\chi^{\vm}\cdot
\kk[\sigma\dual\cap M]
=
\kk[\vpm+(\sigma\dual\cap M)]
$$
by Subsection~(\ref{polToCartier}).
Since we are interested in the degree $m=0$, the upshot is
$$
\gH^0\!\big(U_\sigma,\CL\big)(0)
= \left\{\begin{array}{cl}
\kk\cdot\chi^0 & \mbox{if } 
\,-\vpm = \vm-\vp\in \sigma\dual
\\
0\hspace{1em} & \mbox{otherwise.}
\end{array}\right.
$$

\subsection{A covering for $\Delta^-\setminus\Delta^+$}
\label{DeltaPMSpectral}
To calculate the topological cohomology of $\Delta^-\setminus\Delta^+$,
we will construct a good covering consisting of contractible subsets.
For this
recall that, for each cone $\sigma\in\Sigma$, we have established
the equality $\vp + \sigma\dual=\Delta^++\sigma\dual$
in Subsection (\ref{polNormFan}).
This gives rise to the definition of the open subsets
$$
S(\sigma)\; \coloneqq \; \Delta^-\setminus\big(\vp+\sigma\dual\big) 
\;=\;
            \Delta^-\setminus\big(\Delta^++\sigma\dual\big) 
\;\subseteq\; \Delta^-\setminus\Delta^+.
$$
Now, we are going to prove two claims about these subsets.

\subsubsection{Claim:
$\;\bigcup_{\sigma\in\Sigma}S(\sigma)= \Delta^-\setminus\Delta^+$}
\label{claimUnion}
Let $p\in \Delta^- \setminus \Delta^+$. Then there is an element
$n\in\delta\dual=|\Sigma|$ with $\pairing{p}{n} <
\min\pairing{\Delta^+}{n}$. This $n$ has to be contained
in some cone $\sigma\in\Sigma$, hence
$\min\pairing{\Delta^+}{n} = \pairing{\vp}{\,n}$. 
Hence $p \notin \vp + \sigma\dual$, thus $p\in S(\sigma)$.

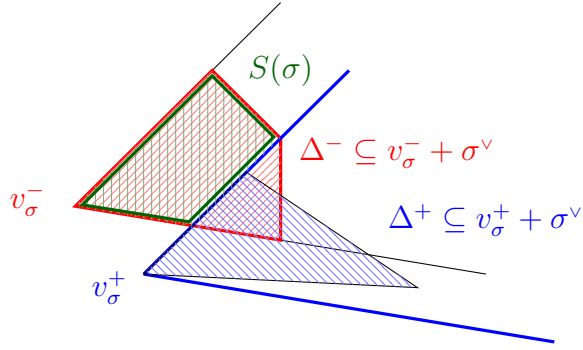
$$
\newcommand{\sizeX}{0.9}
\bild{\sizeX}{
\draw[thin,  color=black]
  (0,0) -- (6,-1) (0,0) -- (3,3);
\draw[very thick,  color=red]  (-0.7,0.1) node{$\vm$};
\draw[very thick, color=red]
  (0,0) -- (3,-0.5) -- (3,1) -- (2,2) -- cycle;
\fill[pattern color=red!50, pattern=north east lines]
  (0,0) -- (3,-0.5) -- (3,1) -- (2,2) -- cycle;
\draw[very thick,  color=red]  (4.7,0.8) node{$\Delta^-\subseteq\vm+\sigma\dual$};
\draw[very thick,  color=blue]
  (1,-1) -- (7,-2) (1,-1) -- (4,2);
\draw[very thick,  color=blue]  (0.5,-1.2) node{$\vp$};
\draw[thin,  color=black]
  (2.5,0.5) -- (5,-1.2) -- (1,-1) -- cycle;
\fill[pattern color=blue!50, pattern=north west lines]
  (2.5,0.5) -- (5,-1.2) -- (1,-1) -- cycle;
\draw[very thick,  color=blue]  (6.0,-0.2) node{$\Delta^+\subseteq \vp+\sigma\dual$};
\fill[pattern color=darkgreen!50, pattern=vertical lines]
  (0,0) -- (1.7,-0.3) -- (3,1) -- (2,2) -- cycle;
\draw[very thick,  color=darkgreen]
  (0.1,0.02) -- (1.66,-0.23) -- (2.90,1.02) -- (2.00,1.92) -- cycle;
\draw[very thick,  color=darkgreen]  (3.0,2.0) node{$S(\sigma)$};
}
$$
\subsubsection{Claim:
If $S(\sigma)\neq\emptyset$, then $\{\vm\}\subseteq S(\sigma)$
is a deformation retract}
\label{claimContract}
If there is some 
$s\in S(\sigma)\subseteq \big(\vm+\sigma\dual\big) \setminus
\big(\vp+\sigma\dual\big)$, 
then there has to be an element $n\in\sigma$ such that
$
\pairing{s}{n}< \pairing{\vp}{n}.
$
Since $\pairing{\vm}{n}\leq\pairing{s}{n}$ 
it follows that $\pairing{s'}{n}< \pairing{\vp}{n}$
for all elements $s'$ of the line segment 
$[{\vm,s}]$.
This implies that $[{\vm,s}]\subseteq S(\sigma)$, 
and altogether this shows that $S(\sigma)$
is star-shaped with respect to $\vm$.

\subsection{Conclusion of the proof}
\label{conclusionEmpty}
By these claims,
$\CS \coloneqq \{S(\sigma)\}_{\sigma\in\Sigma}$
is a \emph{good} covering of $\Delta^-\setminus\Delta^+$ 
in the sense of \cite[\S5]{BottTu}:
First, each set $S(\sigma)$ is open in $\Delta^-\setminus\Delta^+$
and acyclic, i.e.\ $\smash{\tH^*}(S(\sigma))=0$. This holds 
because $S(\sigma)$ is contractible. Second, these properties also hold for
arbitrary intersections, which is obvious here, due to
$S(\sigma) \cap S(\sigma') = S(\sigma\cap\sigma')$.

We build a \v{C}ech complex similarily as in 
Subsection~(\ref{cechSpectral}), but now from relative singular
$0$-chains with coefficients $\kk$
$$
\textstyle
C^p(\CS, \Delta^-\setminus\Delta^+) 
\coloneqq \prod\limits_{i_0<\cdots<i_p}
\gH^0\!\big(\Delta^-\!,\, S(\sigma_{i_0}\cap\ldots\cap\sigma_{i_p})\big).
$$
By \cite[Prop.~10.6]{BottTu}, its $p$-th cohomology is the relative
singular cohomology, which in turn by contractibility of $\Delta^-$,
is isomorphic to the reduced cohomology:
$$
    \gH^p\!\big(C^\bullet(\CS,\Delta^-\setminus\Delta^+)\big) 
= \,\gH^p\!\big(\Delta^-\!,\;\Delta^-\setminus\Delta^+,\kk\big)
= \wt{\gH}^{p-1}\!\big(\Delta^-\setminus\Delta^+\big) .
$$
On the other hand, we obtain from Step \ref{claimContract} that
$$
\gH^0\!\big(\Delta^-\!,\, S(\sigma)\big)=
\wt{\gH}^{-1}\!\big(S(\sigma)\big)= \kk
\;\ifff\;
S(\sigma)=\emptyset
\;\ifff\;
\vm-\vp\in \sigma\dual.
$$
Hence the complexes
$C^\kbb(\CS, \Delta^-\setminus\Delta^+)$
and
$C^\kbb(\CU, \CL)$ from Subsection (\ref{cechSpectral})
coincide, and the claim of Theorem~\ref{th-cohomByDiffB} is proven.

\subsection{Using higher inverse limits}
In Subsection~(\ref{conclusionEmpty}),
we have used \v{C}ech complexes 
mostly for historical reasons -- it was the main tool in all 
proofs of the standard toric cohomology formulas so far.
However, \v{C}ech cohomology is just the down-to-earth description of the
following more abstract, but less technical point of view:
%

The global section functor on $X = \bigcup_{\sigma\in\Sigma} U_\sigma$
is the composition of the local section functors on $U_\sigma$
followed by the inverse limit, i.e.\
$$
\gH^0(X,\CL) = \kprojlimm{\sigma\in\Sigma} \gH^0(U_\sigma,\CL).
$$
Now, the Grothendieck spectral sequence yields
$$
E^{pq}_2(\text{alg}) = \kprojlimm{\sigma\in\Sigma}^p
\gH^q(U_\sigma,\CL) \Longrightarrow \gH^{p+q}(X,\CL)
$$
and, by the same method on the topological side,
$$
E^{pq}_2(\text{top}) = \kprojlimm{\sigma\in\Sigma}^p
 \gH^q(\Delta^-,S(\sigma)) \Longrightarrow 
\gH^{p+q}(\Delta^-,\,\Delta^-\setminus\Delta^+).
$$
Since $E^{pq}_2(\mbox{\rm alg})=E^{pq}_2(\mbox{\rm top})$,
vanishing unless $q=0$, we obtain the result again.

\section{Application: full exceptional sequences of nef line bundles}
\label{app}
\label{sub:fes}

\noindent
Theorem~\ref{th-cohomByDiffB} can be used to quickly check whether a
sequence of nef toric line bundles is exceptional. We illustrate this
on the Hirzebruch surface
$\F_1=\PP\big(\CO_{\PP^1}\oplus\CO_{\PP^1}(1)\big)$
of Subsection~(\ref{easyEx}).

Let $\pi_1 \colon\F_1\to\PP^1$ be the projection as a ruled surface, and 
$\pi_2\colon\F_1\to\PP^2$ be the blowing-down map. The sheaves
$\CLL(a,b) \coloneqq \pi_1^*(\CO_{\PP^1}(a))\otimes_{\CO_{\F_1}}\pi_2^*(\CO_{\PP^2}(b))$ with $a,b\in\Z$ are the toric line bundles on $\F_1$. 
For $a,b\geq 0$, each $\CLL(a,b)$ is nef and corresponds, by 
Subsection~(\ref{polToCartier}), to a lattice polytope $P(a,b)$ in 
$M_\R$.
In the notation of Subsection~(\ref{easyEx}), $A=P(1,0)$ and $B=P(0,1)$,
so that $P(a,b) = aA + bB$.
The ample line bundle $\CLL(1,1)$ corresponds to
  $\Delta \coloneqq P(1,1) = A+B$, i.e.\ $\normal(A+B) = \Sigma$.

\newcommand{\tcp}[2]{\fill[thick] (#1,#2) circle (2pt);} 

\newcommand{\polyt}[6]{%
  \def\sx{#1} \def\sy{#2}
  \draw[color=oliwkowy!40] (-0.3,-0.3) grid (\sx+0.3,\sy+0.3);
  \fill[pattern color=#4!50, pattern=#5] #3 -- cycle;
  \draw[thick,black] #3 -- cycle;
  #6
}

\newcommand{\Paa}{ \polyt{0}{0}{(0,0)}{}{}{\tcp00} }
\newcommand{\Pab}{ \polyt{1}{1}{(0,0)--(0,1)--(1,1)}{black}{north west lines}{}} 
\newcommand{\Pac}{ \polyt{2}{2}{(0,0)--(0,2)--(2,2)}{black}{north west lines}{}} 
\newcommand{\Pba}{ \polyt{1}{0}{(0,0)--(1,0)}{}{}{}{}}                           
\newcommand{\Pbb}{ \polyt{2}{1}{(0,0)--(0,1)--(2,1)--(1,0)}{black}{north east lines}{}} 
\newcommand{\Pbc}{ \polyt{3}{2}{(0,0)--(1,0)--(3,2)--(0,2)}{black}{north east lines}{}} 
\newcommand{\Pcc}{ \polyt{4}{2}{(0,0)--(2,0)--(4,2)--(0,2)}{black}{north east lines}{}}

%
We are looking for a sequence of polytopes
$(P(a_1,b_1),\ldots,P(a_4,b_4))$ corresponding to
a full exceptional sequence of line bundles $(\CL_1,\ldots,\CL_4)$
for the bounded derived category $\Db(\F_1)$,
i.e.\ $\gExt^*(\CL_i,\CL_j)=0$ for $j>i$, and the line bundles generate
$\Db(\F_1)$; see \cite[\S1.4]{Huybrechts}.
%
%
%
%
%
In fact, the following two sequences of polytopes

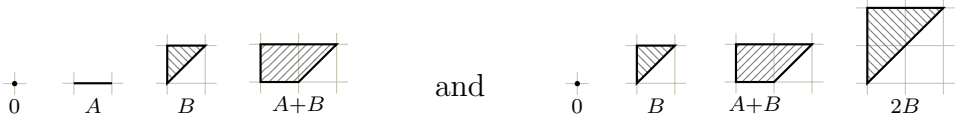
\newcommand{\txtp}[3]{ \node at (#1,#2) {$\scriptstyle #3$}; }

$
\hspace*{\fill}
\bild{0.5}{\Paa \txtp{0.0}{-0.6}{0}}\hspace{1em}
\bild{0.5}{\Pba \txtp{0.5}{-0.6}{A}}\hspace{1em}
\bild{0.5}{\Pab \txtp{0.5}{-0.6}{B}}\hspace{1em}
\bild{0.5}{\Pbb \txtp{1.0}{-0.6}{A+B}}
$
\hfill\raisebox{2ex}{\text{and}}\hfill
$
\bild{0.5}{\Paa \txtp{0.0}{-0.6}{0}}\hspace{1em}
\bild{0.5}{\Pab \txtp{0.5}{-0.6}{B}}\hspace{1em}
\bild{0.5}{\Pbb \txtp{0.5}{-0.6}{A+B}}\hspace{1em}
\bild{0.5}{\Pac \txtp{1.0}{-0.6}{2B}}
\hspace*{\fill}
$

\noindent
give rise to exceptional sequences of line bundles, as is easily
checked with Theorem~\ref{th-cohomByDiffB}.
Indeed, in Subsection~(\ref{easyEx}), we have already seen that
$A$ and $2B$ together cannot be part of an exceptional sequence.
Moreover, these two sequences are full, i.e.\ the respective line bundles
generate $\Db(\F_1)$. This, too, can be seen combinatorially: we use the
general fact that, if $\CO(1)$ is an ample line bundle on a smooth projective
variety $X$ of dimension $d$, then 
$\CO\oplus\CO(1)\oplus\cdots\oplus\CO({d})$
is a strong generator of $\Db(X)$ by \cite[Thm.~4.1]{Orlov-generators}.
{In our example, we need to generate $\CO(2)=\CO(2\Delta)$. We will show 
how to generate} $2\Delta = 2A+2B$ using the
first exceptional sequence, $(0,A,B,A+B)$. This is achieved by the following \emph{polytopal resolution} of $2\Delta$

\newcommand{\rcp}[2]{\fill[thick,red] (#1,#2) circle (4pt);} 

\newcommand{\qolyt}[6]{
  \def\sx{#1} \def\sy{#2}
  \fill[pattern color=#4!50, pattern=#5] #3 -- cycle;
  \draw[thick,black] #3 -- cycle;
  #6
}

\newcommand{\stxtp}[3]{ \node at (#1,#2) {$\smash{\scriptstyle #3}$}; }

\newcommand{\Qaa}{ \qolyt{0}{0}{(0,0)}{}{}{\tcp00 \stxtp{0.0}{-1.0}{0}}}
\newcommand{\Qab}{ \qolyt{1}{1}{(0,0)--(0,1)--(1,1)}{black}{north west lines}
                         {\stxtp{0.2}{-1.0}{B}}}
\newcommand{\Qac}{ \qolyt{2}{2}{(0,0)--(0,2)--(2,2)}{black}{north west lines}
                         {\stxtp{0.5}{-1.0}{2B}}}
\newcommand{\Qba}{ \qolyt{1}{0}{(0,0)--(1,0)}{}{}{}
                         {\stxtp{0.5}{-1.0}{A}}}
\newcommand{\Qbb}{ \qolyt{2}{1}{(0,0)--(0,1)--(2,1)--(1,0)}{black}{north east lines}
                         {\stxtp{0.9}{-1.0}{A+B}}}
\newcommand{\Qbc}{ \qolyt{3}{2}{(0,0)--(1,0)--(3,2)--(0,2)}{black}{north east lines}
                         {\stxtp{1.3}{-1.0}{A+2B}}}
\newcommand{\Qcc}{ \qolyt{4}{2}{(0,0)--(2,0)--(4,2)--(0,2)}{black}{north east lines}{}}

\vspace{-3ex}

\begin{align*}
0 \to
\raisebox{-2.2ex}{%
  \fbox{\bild{0.3}{\Qac \rcp11} }
}
\to
\raisebox{-2.2ex}{%
  \fbox{\bild{0.3}{\Qbc \rcp11}
        \raisebox{0.7ex}{$\oplus$}
        \bild{0.3}{\Qbc \rcp21}
  }
}
\to
\raisebox{-2.2ex}{%
  \fbox{\bild{0.3}{\Qcc \rcp21 \txtp{2.0}{-0.7}{2\Delta} \hspace{1em}} }
}
\to 0
\intertext{where $2B$ and $A+2B$ are in turn resolved by}
\\[-6ex]
0 \to
\raisebox{-0.3ex}{\fbox{\bild{0.3}{\Qba \rcp10}}}
\to
\raisebox{-1.1ex}{%
  \fbox{\bild{0.3}{\Qbb \rcp10}
        \raisebox{0.7ex}{$\oplus$}
        \bild{0.3}{\Qab \rcp11}
  }
}
\to
\raisebox{-1.9ex}{%
  \fbox{\bild{0.3}{\Qac \rcp11}}
}
\to 0
\\
0 \to
\raisebox{-0.5ex}{\fbox{\bild{0.3}{\Qaa \rcp00}}}
\to
\raisebox{-1.3ex}{%
  \fbox{%
    \raisebox{0.8ex}{%
      \raisebox{-0.8ex}{\bild{0.3}{\Qba \rcp10}}
      $\oplus$
      \raisebox{-0.8ex}{\bild{0.3}{\Qba \rcp00}}
      $\oplus$
    }
    \bild{0.3}{\Qab \rcp00}
  }
}  
\to
\raisebox{-1.3ex}{%
  \fbox{\bild{0.3}{\Qbb \rcp00}
        \raisebox{0.8ex}{$\oplus$}
        \bild{0.3}{\Qbb \rcp10}
        \raisebox{0.8ex}{$\oplus$}
        \bild{0.3}{\Qbb \rcp11} }
}
\to
\raisebox{-2.1ex}{%
  \fbox{\bild{0.3}{\Qbc \rcp11} }
}
\to 0
\end{align*}

\noindent
In each polytope, the red circle indicates the origin, thus fixing the position
of the polyhedra in $M_\R$.
Morphisms between polytopes are all possible inclusions, with signs coming
from the natural inclusion-exclusion rule.
The polytopal resolutions induces exact sequences of direct sums of line bundles
by applying the functor
$\PolPlusnullSigma\to\Coh(\F_1)$,
$P \mapsto \CO_X(P)$
of Subsection~(\ref{polToCartier}).

\enlargethispage{1.8cm}

\renewcommand*{\bibliofont}{\footnotesize\rmfamily}
\bibliographystyle{alpha}
\bibliography{dop}

\end{document}